\documentclass[12pt,a4paper]{article}

\usepackage{a4,enumerate}
\usepackage{pstricks,pst-node,pst-coil,pst-plot}
\usepackage{amsthm}
\usepackage{latexsym}

\long\def\regie#1{}

\newdimen\templaenge

\DeclareMathAlphabet{\doba}{U}{msb}{m}{n}

\gdef\mN{\doba{N}}

\gdef\mR{\doba{R}}

\gdef\mZ{\doba{Z}}

\def\qed{{\leavevmode\unskip\nobreak\hfil\penalty 50\hskip 1em%
  \hbox{}\nobreak\hfil\lower 1pt\hbox{$\Box$\kern-.5pt}\parfillskip 0pt
  \finalhyphendemerits 0\par\bigbreak}}
\def\qedmath#1{\setbox0\hbox{$\displaystyle #1$}\templaenge=\textwidth\advance\templaenge by -\wd0%
\setbox1\hbox{$\Box$}\advance\templaenge by -2\wd1%
$$#1\hbox to0pt{\kern.5\templaenge$\Box$\kern-.5pt\hss}$$\par\bigbreak}

\def\al{{\alpha}}
\def\be{{\beta}}

\def\De{{\Delta}}

\def\Om{{\Omega}}
\def\la{{\lambda}}
\def\ka{{\kappa}}

\def\si{{\sigma}}

\def\ga{{\gamma}}
\def\ep{{\varepsilon}}
\def\Ga{{\Gamma}}

\def\ph{{\varphi}}
\def\phi{{\varphi}}

\def\pa{{\partial}}

\def\cK{\mathcal{K}}

\def\cM{\mathcal{M}}

\def\cQ{\mathcal{Q}}

\def\cS{\mathcal{S}}
\def\cW{\mathcal{W}}
\def\cV{\mathcal{V}}

\def\ohne{-}
\def\ti{\tilde}
\def\witi{\widetilde}

\def\ol{\overline}

\def\lan{\langle}
\def\ran{\rangle}

\def\Rplus{\mR^+}

\def\halbplatz{\thinspace}
\let\mo=\mathopen
\let\mc=\mathclose

\def\ie{i.\halbplatz e.\ \ignorespaces}
\def\eg{e.\halbplatz g.\ \ignorespaces}

\def\nummerarray#1#2{\par\noindent\setbox0\hbox{\rm (#1)}\setbox1\hbox{$#2$}\unhcopy0%
\dimen0=.5\textwidth \advance\dimen0 by -\wd0 \advance\dimen0 by -.5\wd1 \kern\dimen0 \unhcopy1}

\def\darea{{\mathop{{\rm darea}}}}
\def\area{{\mathop{{\rm area}}}}

\def\min{\mathop{{\rm min}}}
\def\max{\mathop{{\rm max}}}
\def\osc{\mathop{{\rm osc}}}

\def\torus{{T^2}}

\def\gflach{{g_0}}
\def\gbel{g}
\def\geukl{{g_{\mathop{\rm eucl}}}}
\def\Aflach{{A_0}}
\def\lflach{{l_0}}

\def\res#1#2{{#1}\lower .11ex\hbox{$|$}\lower .644ex\hbox{$\scriptstyle #2$}}

\def\length#1#2{{\mbox{\rm length}_{#2}({#1})}}

\def\cKpo{\cK^{+}}
\def\cKm{\cK^{-}}
\def\sqrr{{\sqrt[\scriptstyle r]{r}\,}}
\def\sys{{\mathop{\rm sys}}_1}
\def\umax{{\max u}}
\def\umin{{\min u}}
\def\uosc{{\osc u}}

\long\def\komment#1{}

\def\proof#1{{\par\medbreak\noindent {\bf Proof\setbox0\hbox{#1}%
\ifdim\wd0=0pt .\else\ \ignorespaces #1.\fi}\enspace}}
\def\iop#1{{\par\medbreak\noindent {\bf Idea of proof\setbox0\hbox{#1}%
\ifdim\wd0=0pt .\else\ \ignorespaces #1.\fi}\enspace}}
\def\examples{{\noindent {\bf Examples. }\par\kern-\baselineskip}}

\newtheoremstyle{remarks}{3pt}{3pt}{}{}{\bfseries}{}{ }{}

\newtheorem{theorem}{\bf T{\footnotesize HEOREM}}[section]
\newtheorem{proposition}[theorem]{\bf P{\footnotesize ROPOSITION}}
\newtheorem{lemma}[theorem]{\bf L{\footnotesize EMMA}}
\newtheorem{corollary}[theorem]{\bf C{\footnotesize OROLLARY}}

\newtheorem*{maintheorem1}{\bf M{\footnotesize AIN} T{\footnotesize HEOREM I}}
\newtheorem*{maintheorem2}{\bf M{\footnotesize AIN} T{\footnotesize HEOREM II}}
\newtheorem{conjecture}[theorem]{\bf C{\footnotesize ONJECTURE}}

\theoremstyle{definition}
\newtheorem*{remark}{Remark}

\newtheorem*{example}{Example}
\newtheorem*{acknowledgement}{Acknowledgement}

\theoremstyle{remarks}

\def\eref#1{{\rm (\ref{#1})}}

\begin{document}
\title{The Willmore Conjecture for immersed tori with small curvature integral}
\author{Bernd Ammann}
\date{June 1999}
\maketitle

\begin{abstract}
\noindent The Willmore conjecture states that any immersion $F:\torus\to \mR^n$
of a 2-torus into euclidean space satisfies $\int_{\torus} H^2\geq 2\pi^2$.
We prove it under the condition that the $L^p$-norm of the 
Gaussian curvature is sufficiently small.
\end{abstract}


{\bf Mathematics Classification:}
53A05 (Primary), 53A30, 58G30 (Secondary)

\section{Introduction}

Let $F:N\to \mR^n$ be an immersion of a closed surface 
$N$ into Euclidean space.
The Willmore integral $\cW(F)$ is defined 
using the mean curvature~$H$ of the immersion~$F$
and the area form~$\darea$ associated 
to the induced metric:
  $$\cW(F):=\int_N H^2\,\darea.$$
In the case $n= 3$ we easily get a lower estimate. If $\ka_1$
and $\ka_2$ are the principal curvatures we have 
$H^2=(\ka_1+\ka_2)^2/4\geq \max\{0,\ka_1\ka_2\}$. 
But $\ka_1\ka_2$ is just the Gaussian curvature $K$ of $(N, F^*\geukl)$ which is
equal to the determinant of the Gauss map $N\to S^2$. Thus we obtain
  $$\cW(F)\geq \int_{N^+}K \, \darea \geq \int_{S^2} 1\,\darea=4\pi,$$
with $N^+=\{x\in N\,|\,K(x)\geq 0\}$ which is mapped onto $S^2$ 
via the Gauss map.

The value $4\pi$ is attained if $F:S^2\to \mR^3$ is the standard embedding. 
And vice versa if $\cW(F)=4\pi$ we know that $N$ is $S^2$ 
whose image $F(S^2)$
is a round sphere.

This well-known result has been improved and generalized by Li and Yau \cite[Fact~3]{li.yau:82}. 
For arbitrary 
dimension $n\geq 3$ we assume that $F^{-1}(p)$ contains $k$ points for some 
$p\in \mR^n$. Then Li and Yau showed the estimate
  $$\cW(F)\geq 4\pi k.$$ 

If $N$ has positive genus, then the value $4\pi$ will never be attained.
In this paper we will study the following conjecture attributed to Willmore \cite{willmore:65}.
 \regie{Kontrollieren!! $n=3$?}
\begin{conjecture}[Willmore conjecture]
For any immersion $F:\torus \to \mR^n$ of the 2-dimensional torus into $\mR^n$, $n\geq 3$, the inequality 
  $$\cW(F)\geq 2\pi^2$$
holds.
\end{conjecture}

Leon Simon proved in \cite{simon:93} that for any fixed dimension $n\geq 3$ the infimum 
  $$\inf\{ \cW(F)\,|\, F:\torus\to \mR^n\}$$
is actually attained and he concludes that 
there is an estimate $W(F)\geq 4\pi+\ep_n$ with $\ep_n>0$ for any $n\geq 3$ without giving 
an explicit value for $\ep_n$.
But the Willmore conjecture remains open until today. 

Nevertheless in many special cases 
the Willmore conjecture has been confirmed. 

If $n=3$ and if the image $F(\torus)$ has a rotational symmetry, 
the Willmore conjecture has been proven by Langer and Singer in
\cite{langer.singer:84}.
Shiohama and Takagi \cite{shiohama.takagi:70} and independently Willmore 
\cite{willmore:71} showed that the Willmore conjecture is true
if $F(\torus)$ is the boundary of an $\ep$-neighborhood of a closed curve in $\mR^3$ with 
$\ep$ sufficiently small.

It should be mentioned here that there is a natural generalization 
of the Willmore functional to immersions $F$ of a closed surface $N$ 
into a Riemannian manifold $(M,h)$ by defining
  $$\cW(F):=\int_N (H^2 + K_M)\,\darea,$$
where $K_M(p)$ is the sectional curvature of $(M,h)$ evaluated at 
the plane $dF(T_pN)$.
This functional is invariant under conformal changes 
of $h$ \cite{thomsen:23,weiner:78} . Hence, the Willmore conjecture for immersions $\torus\to\mR^n$ is equivalent to 
the Willmore conjecture for immersions $F:\torus\to S^n\subset\mR^{n+1}$. The Willmore conjecture for 
immersions of the latter kind
has been proven by Ros \cite{ros:p97} under the additional condition 
that $n=3$ and that $F(\torus)$ is invariant under the 
antipodal map of $S^3$.

\begin{figure}[tbp]
\newgray{liyaufarbe}{.50}
\newgray{montielrosfarbe}{.80}

\def\achseneinst{{\psset{linewidth=1pt,linecolor=gray,linestyle=solid,fillstyle=none}}}
\def\pinpfeileinst{\psset{linewidth=1pt,linecolor=black,linestyle=solid,fillstyle=none}}
\def\modformeinst{\psset{linewidth=2pt,linecolor=black,linestyle=solid,fillstyle=vlines,hatchcolor=black}}
\def\liyaueinst{\psset{linewidth=2pt,linecolor=black,linestyle=solid,fillstyle=vlines*,fillcolor=liyaufarbe,hatchcolor=black}}
\def\montielroseinst{\psset{linewidth=2pt,linecolor=black,linestyle=solid,fillstyle=vlines*,fillcolor=montielrosfarbe,hatchcolor=black}}
\newdimen\descwidth
\descwidth=7.2cm

\begin{center}
\psset{unit=6cm}
\begin{pspicture}(-.2,-.2)(2.2,2.2)
\achseneinst
\psaxes[Dx=.5]{->}(0,0)(-.2,-.2)(1.3,2.0)

\modformeinst
\pscustom{
  \psline(0,1.8)(0,1)
  \psarcn(0,0){1}{90}{60}
  \psline(.5,.866)(.5,1.8)
} 

\liyaueinst
\pscustom{
  \psline(.5,1)(.5,.866)
  \psarc(0,0){1}{60}{90}
  \psline(0,1)(.5,1)
} 

\montielroseinst
\pscustom{
  \psline(.5,1)(.5,1.5)
  \psarc(.5,1){.5}{90}{180}
  \psline(0,1)(.5,1)
} 

\rput(1.42,0){$x$}
\rput(0,2.16){$y$}
\rput[r](.8,.4){\psframe[linecolor=montielrosfarbe,fillstyle=solid,fillcolor=montielrosfarbe](0,-.15)(.3,.15)}
\rput[l](1.2,.4){ 
\vtop{\raggedright\hsize=\descwidth\noindent\small
for tori in these conformal equivalence classes Montiel and Ros proved the Willmore Conjecture}
}
\rput[r](.8,1.1){\psframe[linecolor=liyaufarbe,fillstyle=solid,fillcolor=liyaufarbe](0,-.15)(.3,.15)}
\rput[l](1.2,1.1){ 
\vtop{\raggedright\hsize=\descwidth\noindent\small
for tori in these conformal equivalence classes Li and Yau proved the Willmore Conjecture}
}
\modformeinst
\rput[r](.8,1.8){\psframe(0,-.15)(.3,.15)}
\rput[l](1.2,1.8){
\vtop{\raggedright\hsize=\descwidth\noindent\small
the moduli space of conformal structures}
}
\end{pspicture}
\end{center}
\caption{The moduli space of conformal structures on $\torus$}\label{modulpic}
\end{figure}
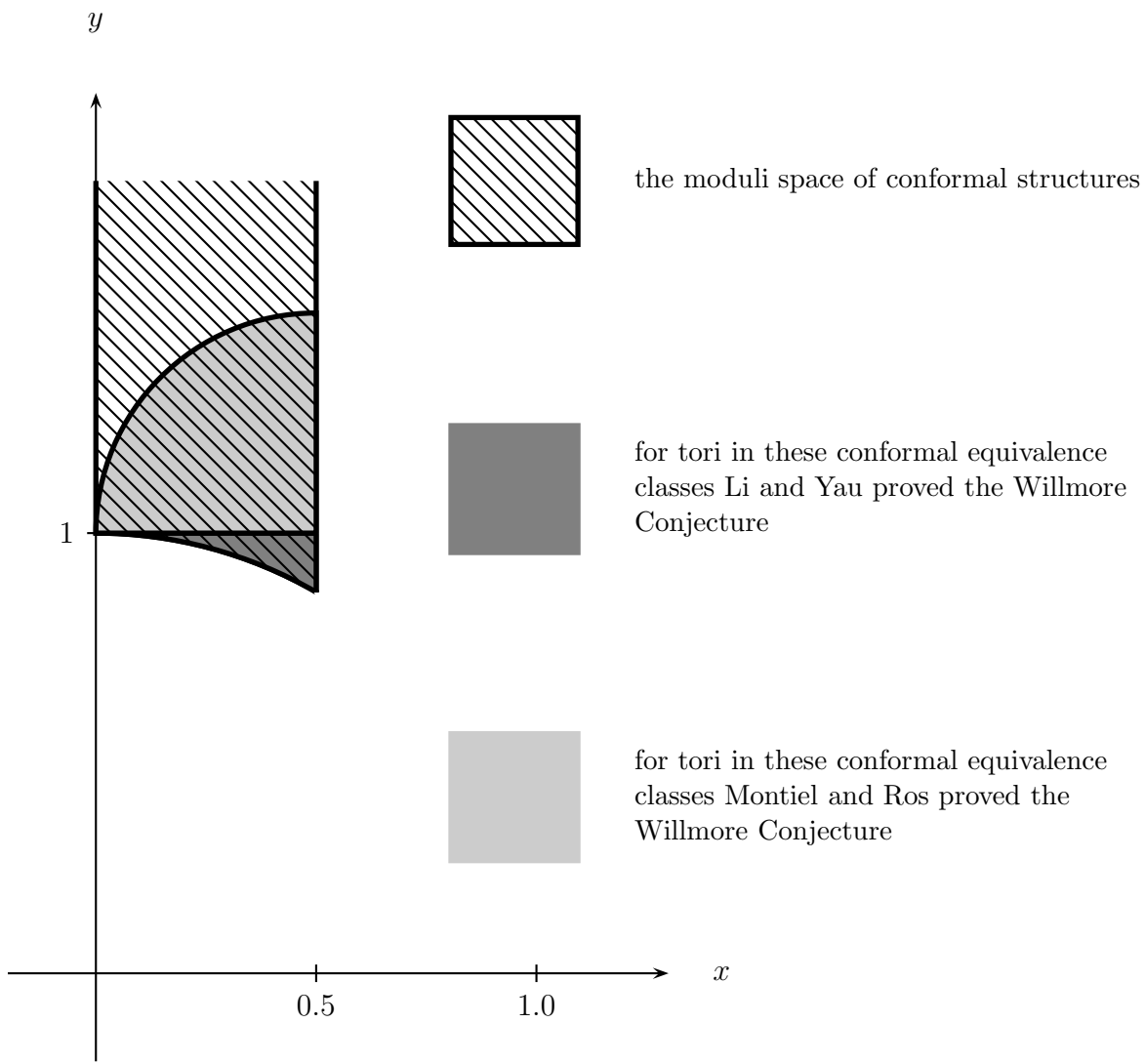

Other partial solutions to the Willmore conjecture use spectral geometry.
If $F:(N,g)\to S^n$ or equivalently $F:(N,g)\to \mR^n$ is a conformal immersion of a closed 
surface $N$ and if $\la_1$ is the first positive eigenvalue of the Laplace operator on $(N,g)$,
then Li and Yau \cite[Theorem~1]{li.yau:82} proved
$$W(F)\geq {1\over 2}\la_1\area(N,g).$$
Every Riemannian 2-torus is conformally equivalent to a flat one, say $(\mR^2/\Ga_{xy},\geukl)$ where the 
lattice $\Ga_{xy}$ is generated by $(1,0)$ and $(x,y)$, $0\leq x \leq 1/2$, $x^2+y^2\geq 1$, $y>0$ and 
where $\geukl$ is the Euclidean standard metric on $\mR^2$.
The first positve eigenvalue of the Laplace operator of $(\mR^2/\Ga_{xy},\geukl)$ is $4\pi^2/y^2$ 
and the area is~$y$.
Thus Li and Yau get the corollary
\begin{equation}\label{liyauhst}
W(F)\geq {2\pi^2\over y}
\end{equation}
which proves the Willmore conjecture for $y\leq 1$, \ie for a subset of the moduli space that has positive measure
(see figure~\ref{modulpic}). The set of conformal equivalence classes for which we know the 
Willmore conjecture has been enlarged by Montiel and Ros 
\cite{montiel.ros:85}. 
They proved that $y\leq 1$ could 
be replaced by the weaker condition 
  $$\left(x-{1\over 2}\right)^2+ (y-1)^2\leq {1\over 4}.$$

Replacing the Laplace operator by the Dirac operator there is a similar result. 
If $F:(N,g)\to S^n$ or $F:(N,g)\to \mR^n$ is an isometric immersion, then  B\"ar \cite{baer:98}
established the estimate
  $$W(F)\geq \mu_1\area(\torus,g)$$
where $\mu_1$ is the first eigenvalue of a twisted Dirac operator on $\torus$.
Unfortunately this statement is no longer true if we replace "isometric immersion" by "conformal immersion", 
so in order to apply this estimate we have to get lower bounds for eigenvalues of 
Dirac operators on non-flat 2-tori. 
Such bounds have been found by the author in \cite{ammanndiss} or alternatively \cite{ammann:p99b}. 
With these estimates we proved the Willmore conjecture 
in an open subset of the spin-conformal moduli space 
provided that the $L^p$-norm of the Gaussian curvature of $N$ is sufficiently small.
This open subset is disjoint from the subset in which Montiel and Ros proved the Willmore conjecture.

A special situation is given if the immersion $F:\torus\to \mR^n$ is flat, \ie the induced metric on 
$F(\torus)$ has vanishing Gaussian curvature. Examples of such tori are the Clifford torus
  $${1\over\sqrt{2}}S^1\times {1\over\sqrt{2}}\subset S^3$$ 
or the Hopf tori \cite{pinkall:85b}. 
For flat immersions the Willmore conjecture is true \cite{chen_bangyen:81}.

Topping \cite{topping:p98} gave a proof of the Willmore conjecture for non-flat immersions
$F:\torus\to S^3\subset\mR^4$ if the following condition is satisfied
  $$\int_{\torus}|K_{\gbel}|\,\darea_\gbel \leq {2\over \pi} \left(W(F)-\area(\torus,\gbel)\right),$$ 
where $\gbel$ is the induced metric. Note that the left hand side is the $L^1$-norm of the 
Gaussian curvature and the right hand side can be rewritten as 
  $${2\over \pi}\,\int_{\torus} \left(H_{\torus\to S^3}\right)^2\darea_\gbel,$$
where $H_{\torus\to S^3}$ is the relative mean curvature of $\torus$ lying in $S^3$.

The main theorems of this article now prove the Willmore conjecture under a similar condition 
on the Gaussian curvature, although our methods are completely different. 
Here we assume that the $L^p$-norm, $p>1$ of the Gaussian curvature 
is bounded by functions that only depend on intrinsic invariants of $(\torus,\gbel)$.
So our assumptions are  --- in contrast to Topping's results ---  purely intrinsic in the following sense:
we construct many non-flat Riemannian metrics $g$ on $\torus$ such that any isometric immersion
$F:(\torus,g)\to(\mR^n,\geukl)$ satisfies $W(F)\geq 2\pi^2$.

\begin{maintheorem1}
For any real number $p>1$ and any conformal equivalence class $c$ on $\torus$ there exists 
$\tau(c,p)>0$ such that the following holds: \newline If $F:\torus \to \mR^n$
is an immersion with induced metric $\gbel:=F^*\geukl$ and conformal 
equivalence class $[\gbel]$ 
satisfying
  $$\left\|K_\gbel\right\|_{L^p(\torus,\gbel)}\cdot \left(\area(\torus,\gbel)\right)^{1-{1\over p}}<\tau([\gbel],p),$$
then the Willmore conjecture
  $$W(F)\geq 2\pi^2$$
holds.
\end{maintheorem1}

\begin{maintheorem2}
There exists a function $\rho:\mo]0,\infty\mc[\times \mo]0,\infty\mc[\times \mo]1,\infty\mc[\to \mo]0,\infty\mc[$ with the following property:
\newline If $F:\torus \to \mR^n$ is an immersion whose 
induced metric $\gbel:=F^*\geukl$ satisfies
  $$\left\|K_\gbel\right\|_{L^p(\torus,\gbel)} < \rho\Bigl(\area(\torus,\gbel),\sys(\torus,\gbel),p\Bigr)$$
for some $p>1$, 
then the Willmore conjecture
  $$W(F)\geq 2\pi^2$$
holds.
\end{maintheorem2}

These results generalize most of the statements about the Willmore conjecture in 
\cite{ammanndiss}. 
The methods used in the proof are strongly related.
The estimate of section~\ref{streckabschsection} can also be used to get 
spectral estimates on $2$-dimensional tori. 
This application will be presented in another article \cite{ammann:p99b} 
which is in preparation. 

\begin{acknowledgement}
I want to thank Christian B\"ar, Ernst Kuwert and Reiner Sch\"atzle 
for many interesting discussions about the subject.
\end{acknowledgement}

\section{Lower bounds for the Willmore functional}
In this section we will prove the main theorems of the article using 
Theorem~\ref{streckabsch} which will be shown in section~\ref{streckabschsection}.

At first we will define some geometric quantities of Riemannian 2-tori $(\torus,g)$.
Let $K_g$ be the Gaussian curvature, $K_g^+:=\max\{K_g,0\}$ and $K_g^-:=\min\{K_g,0\}$. 
For $p\in [1,\infty[$ we set
\begin{eqnarray*}
  \cK_p^{\phantom{\pm}}(g)&:=&
     \left\|K_g^{\phantom{\pm}}\right\|_{L^p(\torus,g)}\cdot \left(\area(\torus,g)\right)^{1-{1\over p}}\\
  \cK_p^\pm(g)&:=&\left\|K_g^\pm\right\|_{L^p(\torus,g)}\cdot \left(\area(\torus,g)\right)^{1-{1\over p}}.
\end{eqnarray*}
The H\"older inequality yields $\cK_p(g)\leq \cK_{p'}(g)$ and 
$\cK_p^\pm(g)\leq \cK_{p'}^\pm(g)$ for $p\leq p'$.
Note that these quantities are invariant under rescaling of the metric.

The \emph{1-systole} or simply the \emph{systole} $\sys(\torus,g)$ is the length of the shortest 
non-contracible loop in $\torus$. We define the geometric quantity
  $$\cV(\torus,\gbel):={\area(\torus,\gbel)\over \sys(\torus,\gbel)^2}$$
which is also invariant under rescaling of the metric. 
Furthermore we set
  $$\cM:=\left\{(x,y)\in\mR^2\,|\,0\leq x\leq 1/2,\;x^2+y^2\geq 1,\;y\geq 0\right\}.$$
As in the introduction we define for $(x,y)\in\cM$ the lattice  
$\Ga_{xy}$ in $\mR^2$ to be generated by $(1,0)$ and $(x,y)$.
Unless otherwise stated $\mR^2/\Ga_{xy}$ always carries
the Riemannian metric induced by the Euclidean metric of $\mR^2$.
Every Riemannian metric on $\torus$ is conformally equivalent 
to exactly one torus of the form 
$\mR^2/\Ga_{xy}$ with $(x,y)\in\cM$. 
Hence $\cM$ can be identified with the moduli space 
of conformal structures on $\torus$. With these notations we have
 $$y=\cV(\mR^2/\Ga_{xy}).$$
The \emph{oscillation} of a continuous function $u:\torus\to\mR$ is defined to be $\osc u:=\max u -\min u$
where the maximum and minimum is to be taken over $\torus$.

Now, we will cite a lemma.
\begin{lemma}[{\cite[1.159]{besse:87}}]\label{konformkruem0lemma}
Let $g_1$ and $g_2=e^{2u}g_1$ be two conformal Riemannian metrics on a surface . 
Then 
their Gaussian curvatures are related via
   $$K_{g_2} - e^{-2u} K_{g_1} = \De_{g_2} u = e^{-2u} \De_{g_1} u.$$ 
\end{lemma}
We always use the convention $\De=-*d*d$, \ie $\De$ has nonnegative 
eigenvalues on compact sets with Dirichlet boundary conditions.

The formula of the lemma also yields a simple proof of the uniformisation theorem 
for $\torus$ stating that any Riemannian
metric on $\torus$ is conformally equivalent to a flat one -- a fact that has 
already been used several times 
in this article.

\begin{lemma}\label{loewnerprep}
Suppose that $\torus$ carries a flat metric $\gflach$ and another metric $\gbel$ conformal to $\gflach$.
Then 
  $$\cV(\torus,\gbel)\geq\cV(\torus,\gflach).$$
\end{lemma}
\proof{}
The proof of this lemma follows a proof of Loewner's theorem, see \eg \cite[~4.1]{gromov:81}.

We write $\gbel=e^{2u}\gflach$. Obviously, $\area(\torus,\gbel)=\int_\torus e^{2u}\,\darea_\gflach$ and
$\area(\torus,\gflach)=\int_\torus \darea_\gflach$.
Let $c:S^1\to \torus$ be a non contractible loop of minimal length $l_0:=\sys(\torus,\gflach)$ 
with respect to $\gflach$. Then for $a\in\torus$ the translated loop 
$c_a(\,\cdot\,):=c(\,\cdot\,)+a$ has the same length with respect to $\gflach$. 
Let $l(a)$ be the length of $c_a$ with respect to
$\gbel$. Then
\begin{eqnarray*}
\int_\torus l(a)\,\darea_\gflach & = & \int_\torus \darea_\gflach \int _{S^1}dt\,\left|\dot{c}_a(t)\right|_\gbel\\
& = & \int_\torus \darea_\gflach \int _{S^1}dt\,e^{u\circ c_a(t)}\,\left|\dot{c}_a(t)\right|_\gflach\\
& = & l_0 \,\int_\torus e^u\,\darea_\gflach\\
& \leq & \l_0\,\area(\torus,\gflach)^{1/2}\,\area(\torus,\gbel)^{1/2}
\end{eqnarray*}
So there is a point $a\in\torus$ with
\qedmath{
  {\sys(\torus,\gbel)^2\over \area(\torus,\gbel)}\leq {l(a)^2\over \area(\torus,\gbel)}
  \leq {{l_0}^2\over \area(\torus,\gflach)}={\sys(\torus,\gflach)^2\over \area(\torus,\gflach)}}
Using the fact that $\cV(\torus,\gflach)\geq \sqrt{3}/2$ for any flat metric $\gflach$ we get the  
\begin{corollary}[Loewner's theorem]
For any Riemannian 2-torus $(\torus,\gbel)$ we have
  $$\cV(\torus,\gbel)\geq {\sqrt{3}\over 2}.$$
\end{corollary}
Equality is attained only for the equilateral flat torus.

\begin{corollary}\label{loewnercor2}
Any isometric immersion $F:(\torus,\gbel)\to(\mR^n,\geukl)$ satisfies
  $$W(F)\geq {2\pi^2\over \cV(\torus,\gbel)}= {2\pi^2\,\sys(\torus,\gbel)^2\over \area(\torus,\gbel)}.$$
In particular, the Willmore conjecture is satisfied for any isometrically immersed torus with 
$\sys(\torus,\gbel)^2\geq \area(\torus,\gbel)$. 
\end{corollary}
\proof{}
We write $\gbel=e^{2u}\gflach$ with $\gflach$ flat. Let $(\torus,\gflach)$ be isometric to 
$\mR^2/\Ga_{xy}$ with $x,y\in\cM$. Then 
  $$y=\cV(\torus,\gflach)\leq \cV(\torus,\gbel).$$ 
Now the corollary follows from inequality~\eref{liyauhst} of the introduction.
\qed

The following lemma is a converse to Lemma~\ref{loewnerprep}.
\begin{lemma}\label{loewnerumk}
Suppose that $\torus$ carries a flat metric $\gflach$ and another metric $\gbel=e^{2u}\gflach$. 
Then 
  $$\cV(\torus,\gbel)\leq e^{2\osc u}\,\cV(\torus,\gflach).$$
\end{lemma}
The proof is straightforward.

We will prove our main theorems in a slightly stronger version than stated 
in the introduction.

\begin{maintheorem1}
There exists a function $\tau:\mo[\sqrt{3}/2,\infty\mc[\times \mo]1,\infty\mc[\to\mo]0,\infty\mc[$ with the 
following property:\newline
If $F:\torus \to \mR^n$ is an immersion such that
the induced metric $\gbel:=F^*\geukl$ satisfies
  $$\cK_p(\gbel)<\tau(y,p)\mbox{\ \ \ and\ \ \  $(T^2,\gbel)$ is conformally equivalent to }\mR^2/\Ga_{xy}$$
for some $p>1$, then the Willmore conjecture
 $$W(F)\geq 2\pi^2$$
holds.
\end{maintheorem1}

\begin{example}
Let $F:(\torus ,\gflach) \to (\mR^n,\geukl)$ be a conformal immersion.
Suppose that ($\torus,\gflach)$ is isometric to $(\mR^2/\Ga_{x2},\geukl)$, $0\leq x\leq 1/2$.
From the explicit construction of $\tau$ in the proof we see that
the Willmore conjecture is satisfied if 
  $$\left\| K\right\|_{L^2(\torus ,\gbel)} 
    \,{\area\Big(F(\torus )\Big)}^{1/2} \leq 0.1987553.$$
\end{example}

\begin{remark}
The $\tau$ constructed in the proof is continuous on $[\sqrt{3}/2,1\mc[\times\mo]1,\infty\mc[$ and
on $\mo]1,\infty\mc[\times\mo]1,\infty\mc[$, but 
  $$\lim_{y\searrow 1} \tau(y,p)=0\neq \tau(1,p)$$
for any $p\in \mo]1,\infty\mc[$. Hence $\tau$ is not continuous at $(y,p)$ with $y=1$.
Nevertheless, if we view $\tau$ as a function on 
$\cM\times \mo]1,\infty\mc[$ we can combine Main Theorem I with the result of 
Montiel and Ros mentioned in the introduction \cite{montiel.ros:85} 
to get a similar function $\witi\tau:\cM\times \mo]1,\infty\mc[\to \mR^+$ 
that is continuous on 
$\left(\cM\ohne\{[\mR^2/\Ga_{01}]\}\right)\times \mo]1,\infty\mc[$ and such
that Main Theorem I holds with $\tau$ replaced by $\witi\tau$. 
\end{remark}

\begin{maintheorem2}
There exists a function 
$\si:\mo[\sqrt{3}/2,\infty\mc[\times \mo]1,\infty\mc[\to \mo]0,\infty\mc[$ 
with the following property:\newline
If $F:\torus \to \mR^n$ is an immersion such that 
the induced metric $\gbel:=F^*\geukl$ satisfies
  $$\cK_p(\gbel)<\si\left(\cV(\torus,\gbel),p\right)$$
for some $p>1$, then the Willmore conjecture 
  $$W(F)\geq 2\pi^2$$
holds.
\end{maintheorem2}

\begin{remark}
In analogy to the previous remark, we cannot chose $\si$ to be continuous at $\cV=1$, $p$ arbitrary. 
But $\si$ can be chosen to be continuous on $\{(\cV,p)\,|\,\cV\neq 1\}$.
\end{remark}

A central role in the proof is played by 
\begin{lemma}\label{beinahflachabsch}
Let $F:\torus\to \mR^n$ be an immersion. Let $\gflach$ denote the standard metric on
$\mR^2$ and suppose that $\torus$ with the induced metric $F^*\geukl$ 
is isometric to $(\mR^2/\Ga_{xy},e^{2u}\gflach)$, $(x,y)\in \cM$, 
where $u$ is a smooth function. Then
  $$W(F)\geq e^{-2 \osc u}\,\pi^2 \Bigl(y+{1\over y}\Bigr)$$
\end{lemma}
This lemma will be shown at the end of this section.

We now define
\begin{eqnarray*}
  \cQ(\cK,p,\cV)& :=  \exp& \Biggl[\;\Big|\log \left(1-{\cK\over 4\pi}\right)\Big|+
  {\cK\over 4\pi -\cK}\,q\log (2q)\\
  &&  \phantom{\Biggl[\;}+ {q\cK\over 2\pi} + {\cK\cV\over 4}\;\Biggr]
\end{eqnarray*}
with $q:=p/(p-1)$.
From the previous lemma and from Theorem~\ref{streckabsch} we get a corollary.

\begin{corollary}
Under the conditions of the previous lemma we have
\begin{eqnarray*}
  W(F)& \geq & \cQ\left(\cK_p(g),p,\cV(\torus,g)\right)^{-1}\,\pi^2 \Bigl(y+{1\over y}\Bigr)\\
  W(F)& \geq & \cQ\left(\cK_p(g),p,y\right)^{-1}\,\pi^2 \Bigl(y+{1\over y}\Bigr).
\end{eqnarray*}
if $\cK_p(g)<4\pi$.
\end{corollary}

\proof{of Main Theorem I}
Let $(\torus,\gbel)$ be conformally equivalent to $\mR^2/\Ga_{xy}$, 
$(x,y)\in\cM$.
We distinguish between two cases. 

In the case $y\leq 1$ we can use the result of Li-Yau \cite{li.yau:82} 
or Montiel-Ros \cite{montiel.ros:85} to see that 
$\tau(y,p)$ can be chosen as any arbitrary positive real number.

On the other hand, for $y>1$ we get 
  $$y+{1\over y}>2.$$
The function $\cQ(\cK,p,\cV)$ is continuous and for
$\cK\to 0$ with $p$ and $\cV$ fixed $\cQ$ converges to~$1$. 
Therefore there is a real number
$\tau(y,p)>0$ such that 
  $$\cQ\left(\cK,p,y\right)^{-1} \left(y+{1\over y}\right)>2$$
whenever $\cK\leq \tau(y,p)$.

Obviously $\tau$ can be chosen as a continuous 
function on $\mo]1,\infty\mc[\times \mo]1,\infty\mc[$.
\qed

\proof{of Main Theorem II}
As in the proof of Main Theorem~I we will distinguish between two cases for the construction of $\si(\cV,p)$.

In the first case, $\cV\leq 1$, according to Corollary~\ref{loewnercor2} the Willmore conjecture is 
satisfied for any immersion $F:\torus\to\mR^n$ with $\cV(\torus,F^*\geukl)=\cV$, 
hence Main Theorem~II
is true if $\si(\cV,p)$ is any positive number.

In the remaining case, \ie $\cV>1,$ we can choose $\si_1(\cV,p)>0$ such that
  $$\cQ(\cK,p,\cV)\leq \sqrt{\cV}$$
whenever $\cK\in [0,\si_1(\cV,p)]$. 
So let $\gbel=e^{2u}\gflach$ be an arbitrary metric on $\torus$ 
with $\gflach$ flat and suppose $\cK_p(\gbel)\leq \si_1(\cV(\torus,\gbel),p)$, then 
according to Theorem~\ref{streckabsch}  
  $$e^{2\osc u}\leq \cQ\Bigl(\cK_p(\gbel),p,\cV(\torus,\gbel)\Bigr)\leq \sqrt{\cV(\torus,\gbel)}.$$  
Using Lemma~\ref{loewnerumk} we get 
  $$\cV(\torus,g_0)\geq \sqrt{\cV(\torus,\gbel)}>1.$$
Let $(\torus,\gbel)$ be conformally equivalent to $\mR^2/\Ga_{xy}$ with $(x,y)\in\cM$. 
Then the above inequality
and Lemma~\ref{loewnerprep} yield 
$\sqrt{\cV(\torus,\gbel)}\leq y\leq \cV(\torus,\gbel)$.

Now set
  $$\si(\cV,p):=\min\biggl\{\si_1(\cV,p),\min\Bigl\{\tau(v,p)\,\Bigm|\,\sqrt{\cV}\leq v\leq \cV\Bigr\}\biggr\}>0$$
with the function~$\tau$ constructed in the proof of Main Theorem~I.
Because of the construction of $\si$  we get
$\si\left(\cV(\torus,\gbel),p\right)\leq \tau(y,p)$.
Therefore Main Theorem II follows from Main Theorem I.
The continuity property is  
clear from the construction of~$\si$.
\qed

\proof{of Lemma~\ref{beinahflachabsch}}
The lemma is a generalization of \cite[Prop.~2, page~287]{li.yau:82}.
We will adapt the proof of this proposition to our situation by introducing 
conformal factors into the formulas.

Let $\gbel:=F^* \geukl$ be the metric on $\torus $ induced 
by the Euclidean metric on $\mR^n$. We decompose $F$ into its coordinate functions
$F=(F_1,\dots,F_n)$, $F_i:\torus \to \mR$. 

We get 
\begin{eqnarray*}
  4 W(F) & = & 4 \int_{\torus } \left|H\right|^2 \,\darea_{\gbel}\\
         & = & \sum_{i=1}^n \int_{\torus }\left(\De_{\gbel}F_i\right)^2\,\darea_{\gbel}\\
         & = & \sum_{i=1}^n \int_{\torus } e^{-2u}\left(\De_{\gflach}F_i\right)^2\,\darea_{\gflach}\\ 
         & \geq &  \sum_{i=1}^n \int_{\torus } e^{-2\umax}\left(\De_{\gflach}F_i\right)^2\,\darea_{\gflach}.
\end{eqnarray*}
After a translation we can assume
  $$\int_{\torus } F_i\,\darea_{\gflach}= 0 \qquad \mbox{for any } i=1,\dots,n.$$
Now we make a Fourier decomposition for the functions $F_i$:
\begin{eqnarray*}
   F_i (w)& = & \sum_{p,q\in\mZ\atop (p,q)\neq (0,0)} A_{ipq}\, \sqrt{2\over y} \,\cos \left( 2\pi 
        \left\lan \left(q, {p-qx\over y}\right),w\right\ran\right)\\
       & & + \sum_{p,q\atop (p,q)\neq (0,0)} B_{ipq} \,\sqrt{2\over y}\, \sin \left( 2\pi 
        \left\lan \left(q, {p-qx\over y}\right),w\right\ran\right).
\end{eqnarray*}
This means for our estimate of $W(F)$
\begin{eqnarray}
  4 \,W(F) & \geq & e^{-2\umax}\,16 \pi^4 \left[\sum_{i,p,q\atop (p,q)\neq (0,0)} 
         A_{ipq}^2\left( q^2+ \left({p-qx\over y}\right)^2\right)^2\right.\nonumber\\ 
  && + \left. \sum_{i,p,q\atop (p,q)\neq (0,0)} 
         B_{ipq}^2\left( q^2+ \left({p-qx\over y}\right)^2\right)^2\right]\nonumber\\
  & \geq & e^{-2\umax} 16 \pi^4 \left[\sum_{i,p,q\atop (p,q)\neq (0,0)} \left(A_{ipq}^2+ B_{ipq}^2\right)
    \left(q^2 + {1\over y^2}\left({p-qx\over y}\right)^2\right)\right].\label{zwifor}
\end{eqnarray}
For the last estimate we used that the inequality
  $$\left[q^2+ \left({p-qx\over y}\right)^2\right]^2 \geq q^2 + \left[2q^2 + 
    \left({p-qx\over y}\right)^2\right]  \left({p-qx\over y}\right)^2 
    \geq q^2 + {1\over y^2}  \left({p-qx\over y}\right)^2
    $$
holds for any $(p,q)\in\mZ\times\mZ$ with $(p,q)\neq (0,0)$ and $x,y\in\mR$ with 
$2y^2\geq 1$.

Now we transform the right hand side of inequality~\eref{zwifor}.
The projections of the vectors $(1,0)$, $(0,1)$ of $\mR^2$ to 
the quotient $\torus =\mR^2/\lan (1,0),(x,y)\ran$ define two vector fields
denoted $e_1$ and $e_2$. These vector fields form an orthonormal frame
for the metric $\gflach$. 

The conformal factor $e^{2u}$ can be calculated as
  $$e^{2u} = \left|\pa_{e_1} F\right|_{\geukl}^2 = \sum_i \left( \pa_{e_1} F_i\right)^2.$$
Therefore we get
  $$\area(\torus ,\gbel)=\int e^{2u}\, \darea_{\gflach} 
                      = \sum_i \int_{\torus } \left(\pa_{e_1}F_i\right)^2\,\darea_{\gflach}
                      = 4 \pi^2 \sum_{i,p,q\atop (p,q)\neq (0,0)} 
                        \left(A_{ipq}^2+ B_{ipq}^2\right) \,q^2.$$
If we replace $e_1$ by $e_2$ we get completely analogously 
  $$\area(\torus ,\gbel)= 4 \pi^2 \sum_{i,p,q\atop (p,q)\neq (0,0)} 
                        \left(A_{ipq}^2+ B_{ipq}^2\right) \,\left({p-qx\over y}\right)^2.$$

On the other hand 
  $$\area(\torus ,\gbel)=\int e^{2u}\,\darea_{\gflach}
    \geq e^{2 \umin}\area(\torus,\gflach)=  e^{2 \umin}\cdot y.$$

All these inequalities together imply
  $$W(F) \geq e^{-2\umax}\, \pi^2\,\left(1+{1\over y^2}\right)
    \area(\torus ,\gbel)\geq e^{-2 \uosc}\, \pi^2 \,\left(y+{1\over y}\right).$$
\qed

\section{Controling the conformal scaling function}\label{streckabschsection}

Consider the 2-dimensional torus with a Riemannian metric $g$.
As in the previous section we write this metric in the form
  $$\gbel= e^{2u} \gflach,$$ 
where $\gflach$ is a flat metric. The function $u$ is uniquely determined by $g$
up to an additive constant.
The aim of this section is to prove estimates for
the oscillation $\osc u = \max u - \min u$.
We will give such an estimate for the case $\cK_1(\gbel)=\int_\torus |K_\gbel|\,\darea_\gbel<4\pi$.
In this situation we can choose $p>1$ such that $\cK_p(\gbel)<4\pi$.
In order to motivate the condition $\cK_p(\gbel)<4\pi$ we construct some examples in 
section~\ref{exsec}.
If $\cK_p(\gbel)$ is sufficiently big, then there is a bubbling phenomenon 
as described for example in \cite{chen_xiuxiong:98} 
for the special case $p=2$. 

\begin{theorem}\label{streckabsch}
We assume for some $p\in \;\mo]1,\infty\mc[$ 
  $$\cK_p:=\cK_p(\torus ,\gbel)<4\pi.$$
Then the oscillation of $u$ is bounded as follows
\nummerarray{a}{\osc u \leq \cS\Bigl(\cK_p,p,\cV(\torus,\gflach)\Bigr),}
\nummerarray{b}{\osc u \leq \cS\Bigl(\cK_p,p,\cV(\torus,\gbel)\Bigr),}
\par\noindent
where we use the definition
  $$\cS(\cK,p,\cV):={1\over 2}\bigg|\log \left(1-{\cK\over 4\pi}\right)\bigg|+
    {\cK\over 8\pi -2\cK}\,q\log(2q)
    + {q\cK\over 4\pi} + {\cK\cV\over 8}$$
with $q:=p/(p-1)$.  
\end{theorem}
\label{schrankedef}
 
\begin{remark}
  $$\cQ(\cK,p,\cV) = \exp \Bigl(2\cS(\cK,p,\cV)\Bigr).$$
\end{remark}

\begin{remark}
The bounds we get are related to a result of Brezis and Merle \cite{brezis.merle:91}
proving the existence of a bound for $\|u\|_{L^q(\mR^2,\gflach)}$ for functions $u:\Om\to\mR$ defined
on a bounded domain $\Om\subset \mR^n$ satisfying the Kazdan-Warner equation 
  $$\De_\gflach u = K(x)e^{2u}$$ 
with Dirichlet boundary conditions.
Their bound depends on $\|K_\gbel\|_{L^p(\mR^2,\gflach)}$ and $\Om$. 
The main difference to our results is that Brezis and Merle do not treat functions $u:\torus\to\mR$
and that their norms are taken with respect to the flat metric $\gflach$ whereas we bound $\osc u$
in terms of the $L^p$-norm with respect to the metric $\gbel=e^{2u}\gflach$.  
\end{remark}

Now we turn to the proof of Theorem~\ref{streckabsch}. The theorem 
will follow from some propositions and corollaries that will take the 
rest of this section. For the proof we split the torus $\torus$
into three subsets two of which have the property that their 
fundamental group is mapped trivially to $\pi_1(\torus)$.
On these two subsets $\osc u$ can be estimated using Corollary~\ref{uupperboundcor}
and Proposition~\ref{ulowerbound}. The third part contains generators of $\pi_1(\torus)$ and 
will be dealt in Proposition~\ref{umidbound}.

\begin{proposition}\label{uupperbound}
Let $G$ be a bounded open subset of $\mR^2$ with the standard 
metric~$\gflach$. 
For a smooth function $u:\ol{G}\to \mR$ with $u\geq 0$, $\res{u}{\partial G}\equiv 0$
we set $\gbel=e^{2u}\gflach$.
We define $\mu_0:=\area(G,\gbel)$ and for the Gaussian curvature  
$K_\gbel$ with respect to $\gbel$
let  
  $$k(A):=\sup \left\{\int_{\widehat G, \gbel} K_\gbel\; 
    \Biggm|\;\widehat G \mbox{\rm \  open subset of $G$,\ } \area(\widehat G,\gbel)=A\right\}.$$
If there are $\ka\in\;]0,2\pi[$, $C>0$, $r\in\;]0,1]$ and $\mu_1\in\;]0,\mu_0]$, such that
$$k(A)\leq \left\{\matrix{C\cdot A^r \leq \ka & \mbox{for}&  0\leq A \leq \mu_1\hfill\cr
                          \ka \hfill &\mbox{for} &  \mu_1 \leq A \leq \mu_0, \hfill }\right.$$
then
  $$\max u \leq {1\over 2}\left|\log \left(1-{\ka\over 2\pi}\right)\right|+ 
    {\ka\over 4\pi -2\ka}\log \left({\mu_0\over \sqrr\mu_1}\right).$$
\end{proposition}

Before proving this proposition, we will prove a corollary.

\begin{corollary}\label{uupperboundcor}
Suppose the open subset $G$ of $\torus$ has the property that any loop $\ga:S^1\to \overline{G}$  
is contractible in $\torus$. 
Let $\gflach$ be a flat metric on $\torus$ 
and $\gbel=e^{2u}\gflach$ another metric on $\torus$ with a smooth function 
$u:\torus\to \mR$ satisfying $\res{u}{\partial G}\equiv 0$ and 
$\res{u}{G}\geq 0$.
Suppose for some $p\in\;\mo]1,\infty\mc[$ we have
  $$\cK_p=\cK_p(\torus ,\gbel)<4\pi.$$
Then we get the estimate
  $$\max_{x\in\ol G} u(x) \leq {1\over 2}\,\biggl|\log \left(1-{\cK_p\over 4\pi}\right)\biggr|\,+\,
    {\cK_p\over 8\pi -2\cK_p}\,q\log (2q),$$
with $q:=p/(p-1)$.
\end{corollary}

\proof{of Corollary~\ref{uupperboundcor}}
Because of the Gauss-Bonnet theorem we have
$$ \left|\int_{\widehat G,\gbel}K_\gbel \right|= \left|\int_{\torus \ohne\widehat G,\gbel}K_\gbel \right|$$
for any open subset $\widehat G\subset G$. 
Therefore
$$ \left|\int_{\widehat G,\gbel}K_\gbel \right|\leq {1\over 2}  \int_{\torus ,\gbel}\left|K_\gbel \right|
   \leq {1\over 2}\,\cK_p.$$
On the other hand, if we have $\area(\widehat G,\gbel)<2^{-q}\area(\torus,\gbel)$, then the estimate
$$ \left|\int_{\widehat G,\gbel}K_\gbel \right|\leq 
   \left\|K_g\right\|_{L^p(\widehat G,\gbel)}
   \left(\area(\widehat G,\gbel)\right)^{1/q}\leq 
   \cK_p\left({\area(\widehat G,\gbel)\over \area(\torus ,\gbel)}\right)^{1/q}$$
is better.

Since any loop $\ga:S^1\to \overline{G}$ is contractible in $\torus$, we can lift $G$ to 
$\mR^2$.

We can apply Proposition~\ref{uupperbound} with
$\ka:=(1/2)\,\cK_p$, $r:=1/q$, $C:=\cK_p\,\area(\torus ,\gbel)^{-1/q}$, 
$\mu_0:=\area(G,\gbel)\leq \area(\torus ,\gbel)$ and 
$\mu_1:=\min\{\mu_0,2^{-q}\area(\torus ,\gbel)\}$.
\qed

\proof{of Proposition~\ref{uupperbound}}
For $v\in \mR$ we define $G(v):=\{x\in G\,|\,u(x)>v\}$. 
The area of $G(v)$ with respect to $\gbel$ (or $\gflach$ resp.)
will be denoted $A(v)$ (or $\Aflach(v)$ resp.). 
For the length of the boundary $\pa G(v)$ 
we write $l(v)$ (or $\lflach(v)$ resp.). 
The functions $A(v)$, $\Aflach(v)$, $l(v)$ and $\lflach(v)$
are differentiable at every regular value $v$ of the function $u$.
Therefore we get for regular values $v$:
\begin{eqnarray}
  e^{2v}\Aflach(v) & \leq &  A(v)\; \leq\; e^{2\umax} \Aflach(v) \label{avunglei}\\
  e^{2v}{d\over dv} \Aflach(v) & = & {d\over dv}A(v)\label{davglei}\\ 
  e^v\lflach(v) & = & l(v) \label{lvglei}\\
  \int_{G(v),\gbel}K_\gbel & \leq  & k(A(v)).
  \label{kint1}
\end{eqnarray}
On the other hand, according to Lemma~\ref{konformkruem0lemma} we have
\begin{equation}
\int_{G(v),\gbel}K_\gbel = \int_{G(v),\gbel}\De_{\gbel} u = - \int_{\pa G(v)} *\, du = 
    \int_{\pa G(v), \gflach} |du|_\gflach.\label{kint2}
\end{equation}
The last equation follows since $u$ is equal to $v$ on $\pa G(v)$ and greater than $v$ on $G(v)$ 
and therefore 
$*\,du$ is negatively oriented on $\pa G(v)$. We also see that $\int_{G(v),\gbel}K_\gbel$ 
is positive for every regular value $v\in u(G)$ of the function $u$.

Using \eref{kint1} and \eref{kint2} we calculate
\begin{eqnarray}
   -{d\over dv} \Aflach(v) & = & \int_{\pa G(v),\gflach} {1\over |du|_\gflach} 
   \;\geq\; {{\lflach(v)}^2\over  \int_{\pa G(v), \gflach} |du|_\gflach} \nonumber \\
   & \geq & {{\lflach(v)}^2\over k(A(v))}\label{davabsch}
\end{eqnarray}  
for any regular value $v\in u(G)$ of $u$.
Now we apply the isoperimetric inequality
\begin{equation}\label{isoper}
{\lflach(v)}^2\geq 4 \pi \, \Aflach(v)
\end{equation}
and we get 
\begin{equation}
  -{d\over dv} \Aflach(v)\geq {4 \pi\over k(A(v))}\Aflach(v). \label{vorfu}
\end{equation}

We set $\mu_2:=\sqrr\mu_1$ and $u_2:=\inf\Bigl\{\al\in[0,\max u\mc[\;\Big|\;A(\al)\leq \mu_2\Bigr\}$.
Let us distinguish between the two cases $v\geq  u_2$ and $v<u_2$. 

We start with $v\geq u_2$. In this case the inequalities 
\eref{avunglei}, \eref{davglei}, \eref{vorfu}
and the bound of $k(A)$ provide the estimate
  $$-{d\over dv}{A(v)}\geq {4\pi\over C\cdot {A(v)}^r} \,e^{2(v-\umax)} A(v),$$
and therefore
$$
   -{1\over r}{d\over dv}\left({A(v)}^r\right)\geq {4\pi \over C}\,e^{2(v-\umax)}.
$$
We use the following lemma. 
Note that any monotonically increasing function is differentiable almost everywhere.

\begin{lemma}\label{strealemma}
Let $f,g:[a,b]\to \mR$ be functions with $f$ monotonically increasing, $g$ 
continuously differentiable and $f(a)=g(a)$. If $f'\geq g'$ almost everywhere, then $f(b)\geq g(b)$.  
\end{lemma} 

This lemma ensures that we can integrate the previous inequality from $u_2$ to $\umax$. 
We use $A(\umax)=0$ and get
$$ {1\over r}{A(u_2)}^r \geq {2\pi \over  C}\, e^{-2\umax}\left(e^{2\umax} -e^{2u_2}\right)
   = {2\pi \over  C}\left(1 -e^{2(u_2-\umax)}\right).$$
Furthermore $A(u_2)=\lim\limits_{\al\searrow u_2} A(\al)\leq \mu_2=\sqrr\mu_1$ and $C{\mu_1}^r\leq \ka$ yield
$$ 1-e^{2(u_2-\umax)} \leq {C\over 2\pi r}{A(u_2)}^r \leq {\ka\over 2\pi}$$
and therefore
\begin{equation}
(\umax) - u_2 \leq {1\over 2}\,\biggl|\log \left(1-{\ka\over 2\pi}\right)\biggr|.\label{umaxschr}
\end{equation}

Now we treat the case $v<u_2$. 
From the estimate~\eref{umaxschr} we know for small $\ep>0$
  $$\Aflach(u_2-\ep) \geq e^{-2\umax}A(u_2-\ep)\geq e^{-2\umax}\mu_2\geq e^{-2u_2}\mu_3,$$ 
with
  $$\mu_3:=\left(1-{\ka\over 2\pi}\right)\mu_2.$$
Using \eref{davglei} we obtain 
  $$\Aflach(v)\geq \int_v^{u_2-\ep} e^{-2\bar v}\left(-{d\over d\bar v}{A(\bar v)}\right)d\bar v \, 
    + \,e^{-2u_2} \mu_3.$$
Here we used the fact that any monotonically decreasing function $h:[a,b]\to \mR$ satisfies
$h(a)-h(b)\geq \int_a^b -h'(t)\, dt$.
In particular we can integrate over the singularities.

Inequality~\eref{vorfu} then provides 
\begin{eqnarray*}
-{d\over dv}{A(v)} & \geq & {4\pi\over \ka}\, e^{2v}
    \left(\int_v^{u_2-\ep} e^{-2\bar v}\left(-{d\over d\bar v}{A(\bar v)}\right)d\bar v \, 
    + \,e^{-2u_2} \mu_3\right).
\end{eqnarray*}
Let $f$ be the solution of the integral equality corresponding to this integral inequality, 
\ie 
\begin{eqnarray*}
  f(v)& = & {4\pi \over  \ka}\, e^{2v} 
    \left(\int_v^{u_2-\ep} e^{-2\bar v}f(\bar v)\, d\bar v \, + \,e^{-2u_2} \mu_3\right).
\end{eqnarray*}
Via differentiation we get the differential equation
  $${d\over dv} f (v) = 2 f(v) -  {4\pi \over  \ka}f(v)$$
with initial value 
$f(u_2-\ep)=(4\pi/\ka)e^{-2\ep}\mu_3$. 
So the solution is
  $$f(v)= {4\pi\over \ka}\,e^{-2\ep}\,\mu_3 \,e^{\left({4\pi\over \ka}-2\right)\left(u_2-v-\ep\right)}.$$
Here we use another elementary lemma. 
\begin{lemma}\label{strealemma2}
Let $f_1$ and $f_2$ be $L^1$-functions on $[a,b]$ and 
$g_1,g_2:[a,b]\to \Rplus$ continuous functions. Let $C\in \Rplus$.
We assume that for any $t\in [a,b]$ we have
\begin{eqnarray*}
  f_1(t) & \geq & g_1(t)\left(C+\int_a^t g_2(s) f_1(s) \,ds\right)\\
  f_2(t) & =    & g_1(t)\left(C+\int_a^t g_2(s) f_2(s) \,ds\right).\\  
\end{eqnarray*}
Then we get $f_1(t)\geq f_2(t)$ for any $t\in[a,b]$. 
\end{lemma}

Thus we obtain
\begin{eqnarray*}
-{d\over dv}{A(v)} & \geq & f(v) = {4\pi\over \ka}\,\mu_3\,e^{-2\ep}\,e^{\left({4\pi\over \ka}-2\right)\left(u_2-v-\ep\right)}\\
\end{eqnarray*}
and the limit $\ep\to 0$ yields 
\begin{eqnarray*}
-{d\over dv}{A(v)} & \geq & f(v) = {4\pi\over \ka}\,\mu_3\,e^{\left({4\pi\over \ka}-2\right)\left(u_2-v\right)}\\
\end{eqnarray*}
Now integration from $0$ to $u_2$ using Lemma~\ref{strealemma} yields
  $$\mu_0-\mu_2 \geq A(0)-A(u_2) \geq \mu_3{\left(1- {\ka \over 2\pi}\right)}^{-1}\left(e^{\left({4\pi\over \ka}-2\right)u_2}-1\right).$$ 
Hence
\begin{eqnarray*}
  u_2& \leq & {\ka\over 4\pi - 2\ka}\log \left(1+ \left(1- {\ka \over 2\pi}\right) {(\mu_0-\mu_2)\over\mu_3}\right)\\
     &   =  & {\ka\over 4\pi - 2\ka}\log \left({\mu_0\over\sqrr\mu_1}\right).
\end{eqnarray*}
This together with \eref{umaxschr} yields the estimate of the proposition. 
\qed

\begin{proposition}\label{ulowerbound}
Let $u:\ol{G}\to \mR$, $u\leq 0$ be a smooth function on the bounded open subset $G\subset \mR^2$
satisfying the boundary condition $\res{u}{\partial G}\equiv 0$. 
Let $\gflach$ be the restriction of the standard metric on $\mR^2$.
We set $\gbel:=e^{2u}\gflach$.
For $p\in\;\mo]1,\infty\mc[$ we define $q:=p/(p-1)$.
Then
  $$\min u \geq - {q\,\cKm_p(G,\gbel)\over 4\pi}.$$
\end{proposition}

\proof{of Proposition~\ref{ulowerbound}}
This time we define $\ka:=\cKm_p(G,\gbel)$
and $G(v):=\{x\in G\,|\,u(x)<v\}$. As in the proof of Proposition~\ref{uupperbound} let
$A(v)$ and $l(v)$ or $\Aflach(v)$ and $\lflach(v)$ resp.\ 
be the area of $G(v)$ and 
the length of $\pa G(v)$ with respect to $\gbel$ or $\gflach$ resp.
Again we have \eref{davglei}, \eref{lvglei} and \eref{isoper}
whereas we have to modify \eref{avunglei} and \eref{kint2}:
\begin{eqnarray}
  e^{2v}\Aflach(v) & \geq & A(v)\label{avunglei2}\\
 \int_{G(v),\gbel}K_\gbel & = & - \int_{\pa G(v), \gflach} |du|_\gflach.\nonumber   
\end{eqnarray}
Furthermore we get
\begin{equation}\label{knabsch}
  k(v):=-\int_{G(v),\gbel}K_\gbel
    \leq \left\|K^-_\gbel\right\|_{L^p(G(v)),\gbel}\area(G(v),\gbel)^{1/q}
    \leq \ka \left({A(v)\over A(0)}\right)^{1/q}.
\end{equation}

Inequality~\eref{davabsch} holds with a different sign:
\begin{equation}
{d\over dv}\Aflach(v)\geq {{\lflach(v)}^2 \over k(v) }.\label{vdavm}
\end{equation}
Together with \eref{isoper} and \eref{knabsch} this yields
$$
   {d\over dv} \Aflach(v) 
   \geq  {4\pi\over \ka}\,\Aflach(v)\left(A(0)\over A(v)\right)^{1/q}.
$$
Finally with \eref{davglei} and \eref{avunglei2} we obtain
$$
  {d\over dv} A(v) \geq {4\pi\over \ka} A(v) \left(A(0)\over A(v)\right)^{1/q}.$$ 
We use again Lemma~\ref{strealemma} in order to integrate 
this inequality from $\umin$ to $0$,
and we get
  $$q\left(A(0)^{1/q}-A(\umin)^{1/q}\right)\geq {4\pi\over \ka} A(0)^{1/q}\,|\umin|.$$
Since $A(\umin)=0$ this implies
\qedmath{|\umin| \leq {q\ka\over 4\pi}.}

\begin{lemma}\label{boundest}
Let $G$ be an open set in $(\torus ,\gbel)$ with smooth boundary. 
We suppose that there are closed curves 
$c_1:[0,1]\to G$ and $c_2:[0,1]\to \torus \ohne \ol{G}$, whose 
corresponding homology classes in $H_1(\torus ,\mZ)$ are not zero. 
Then
  $$\length{\pa G}{\gbel}\geq 2 \,\sys(\torus ,\gbel).$$ 
\end{lemma}

\proof{}
If we regard $\overline{G}$ as a $2$-cycle, then clearly $\pa G$ is homologous to zero.
Now decompose $\pa G$ into its components. Each component is 
diffeomorphic to $S^1$. 

We will show that there is at least one component non-homologous to zero. Together 
with $[\pa G]=0$ this implies that are at least two such components and therefore we get
the statement of the lemma.

So let us suppose that all components of $\pa G$ are homologous to zero.
Let $\pi:\mR^2\to \torus $ be the universal covering. 
Then $\pi^{-1}(\pa G)$ is diffeomorphic to a disjoint union of countably many $S^1$.
We write
  $$\pi^{-1}(\pa G)=\dot{\bigcup_{i\in \mN}} Y_i$$ 
with $Y_i\cong S^1$.
We choose lifts $\ti c_i:\mR\to \mR^2$ of $c_i$, \ie
$\pi\left(\ti c_i(t+z)\right)=c_i(t)$ for all $t\in [0,1]$, $z\in \mZ$ and $i=1,2$.
Then we take a path $\ti\ga:[0,1]\to \mR^2$ joining $\ti c_1(0)$ to $\ti c_2(0)$.
We define $I$ to be the set of all $i\in \mN$ such that $Y_i$ meets the trace of $\ti\ga$. 
The set $I$ is finite.
Using the Theorem of Jordan and Schoenfliess about simply closed curves in $\mR^2$ we can  
inductively construct a compact set $K\subset \mR^2$ with boundary $\bigcup_{i\in I}Y_i$. 
Either $\ti c_1(0)$ or $\ti c_2(0)$ is in the interior of $K$.
But if $\ti c_i(0)$ is in the interior of $K$, then the whole trace  $\ti c_i(\mR)$ 
is contained in $K$. Furthermore,  
$\ti c_i(\mR)=\pi^{-1}\left(c_i([0,1])\right)$ is closed and therefore
compact. This implies that $c_i$ is homologous to zero in contradiction to our assumption.
\qed

\begin{proposition}\label{umidbound}
Suppose that $\torus $ carries a Riemannian metric $\gbel = e^{2u}\gflach$ with $\gflach$ 
flat. 
Let $c_i:S^1\to \torus $, $i=1,2$ be closed paths non-homologous to zero.
We set
  $$v_1:=\max_{t\in S^1} u\circ c_1(t) \quad \mbox{and} \quad v_2:=\min_{t\in S^1} u\circ c_2(t).$$
Then 
\nummerarray{a}{v_2-v_1 \, \leq \,{\displaystyle\cK_1(\torus,\gbel)\, \cV(\torus,\gflach)\over \displaystyle \lower1mm\hbox{$8$}},}
\nummerarray{b}{v_2-v_1 \, \leq \,{\displaystyle\cK_1(\torus,\gbel)\, \cV(\torus,\gbel)\over \displaystyle \lower1mm\hbox{$8$}}.}
\end{proposition}

\proof{}
The statement is void for $v_2\leq v_1$, therefore we can assume $v_2>v_1$.
We set $\cK_1:=\cK_1(\torus ,\gbel)$.
Let $v\in\;]v_1,v_2[$ be a regular value of $u$.
As in the proof of the previous proposition we set $G(v):=\{x\in\torus \,|\,u(x)<v\}$, 
let $A(v)$ be the area $G(v)$ with respect to $\gbel$, 
$l(v)$ the length of $\pa G(v)$ with respect to $\gbel$, 
and $\Aflach(v)$ and $\lflach(v)$ the area and length with respect to $\gflach$.
In analogy to \eref{kint2} we get 
  $$ \int_{G(v),\gbel}K_\gbel =  - \int_{\torus  \ohne G(v),\gbel}K_\gbel = 
      -\int_{\pa G(v)} *\, du  $$
and therefore
  $$\int_{\pa G(v), \gflach} |du|_\gflach=\int_{\pa G(v), \gbel} |du|_\gbel=\int_{\pa G(v)} *\, du
    \leq {1\over 2}\int_{\torus ,\gbel}\left|K_\gbel\right|\leq {\cK_1\over 2}.$$
We obtain
  $$ {d\over dv} \Aflach(v)  =  \int_{\pa G(v),\gflach} {1\over |du|_\gflach} 
   \geq {{\lflach(v)}^2\over  \int_{\pa G(v), \gflach} |du|_\gflach}\geq 2\,{{\lflach(v)}^2\over \cK_1}.$$
Using Lemma~\ref{boundest} we know that the right hand side of this inequality  
is greater than or equal to
$8\,{\sys(\torus,\gflach)}^2/\cK_1$.

Integration yields 
  $$\area(\torus ,\gflach)\geq \Aflach(v_2)-\Aflach(v_1)\geq 8 \,{{\sys(\torus,\gflach)}^2\over \cK_1}\,(v_2-v_1),$$
and therefore the statement of~(a). 

Similarly, we show statement~(b):
  $$ {d\over dv} A(v)  =  \int_{\pa G(v),\gbel} {1\over |du|_\gbel} 
   \geq {{l(v)}^2\over  \int_{\pa G(v), \gbel} |du|_\gbel}\geq 8 \,{{\sys(\torus,\gbel)}^2\over \cK_1},$$
\qedmath{\area(\torus ,\gbel)\geq A(v_2)-A(v_1)\geq 8 \,{{\sys(\torus,\gbel)}^2\over \cK_1}\,(v_2-v_1).}
An alternative way to prove (b) is to use (a) together with Lemma~\ref{loewnerprep}.

\proof{of Theorem~\ref{streckabsch}}
Again we set $G(v):=\{x\in\torus \,|\,u(x)>v\}$. 
Let~$v_2$ be the supremum of all $v\in \mR$ with the property that there is a closed path 
$c_2(v):S^1\to G(v)$
that is non-homologous to zero in $\torus$.
Similarly, we define $\widehat G(v):=\{x\in\torus \,|\,u(x)<v\}$, 
and let $v_1$ be the infimum of all $v\in \mR$ for which there is a closed path
$c_1(v):S^1\to\widehat G(v)$ that is non-homologous to zero in $\torus$.

For any $\ep>0$ we have $\max \bigl(u\circ c_1(v_1+\ep)\bigr) <v_1+\ep$
and $\min \bigl(u\circ c_2(v_2-\ep)\bigr)> v_2-\ep$. 
We apply Proposition~\ref{umidbound}, and the limes
$\ep\to 0$ yields
\begin{eqnarray}
  v_2-v_1 & \leq &{\cK_1(\torus ,\gbel)\, \cV(\torus,\gflach)\over 8}\nonumber\\
  &\leq & {\cK_p(\torus ,\gbel)\, \cV(\torus,\gflach)\over 8}.
   \label{streasummand1}
\end{eqnarray}
Now we apply Corollary~\ref{uupperboundcor} for $G:=G(v)$ and $v:=v_2+\ep$ where we 
replace $u$ by $u-v$.
We get in the limit $\ep\to 0$
\begin{equation}
(\umax)-v_2\leq {1\over 2}\biggl|\log \left(1-{\cK_p\over 4\pi}\right)\biggr|+
    {\cK_p\over 8\pi -2\cK_p}\,q\log (2q),
     \label{streasummand2}
\end{equation}
with $q:=p/(p-1)$.

Similarly, Proposition~\ref{ulowerbound} yields
\begin{equation}
v_1-(\umin)\leq {q\,\cK_p(G,\gbel)\over 4\pi}.\label{streasummand3}
\end{equation}

Adding inequalities~\eref{streasummand1}, \eref{streasummand2} and \eref{streasummand3} 
we obtain statement~(a).

The proof of statement~(b) is completely analogous.
\qed

\section{An estimate on the disk}

From the results of the previous section is not difficult to derive another theorem. 

\begin{theorem}\label{disktheo}
Let $\gbel$ be a Riemannian metric on a domain $G$ whose closure 
is diffeomorphic to the $2$-dimensional disk. 
We write $\gbel$ as $\gbel=e^{2u}\gflach$ with $\gflach$ flat and 
$\res{u}{\partial G}\equiv 0$.
For $p\in\;\mo]1,\infty\mc[$ we assume $\cKpo_p:=\cKpo_p(G,\gbel)<2\pi$.
Then we get the estimate
  $$\max u \leq {1\over 2}\,\Biggl|\log \left(1-{\cKpo_p\over 2\pi}\right)\Biggr|\,+\,
    {\cKpo_p\over 4\pi -2\cKpo_p}\,q\log q,$$
with $q:=p/(p-1)$.
\end{theorem}
 
\proof{of Theorem~\ref{disktheo}}
W.l.o.g.\ we can assume $u\geq 0$.

For any open subset $\widehat{G}\subset G$ we have
\begin{eqnarray*}
\int_{\widehat G,\gbel}K_\gbel & \leq  & \left\|K^+_\gbel\right\|_{L^1(\widehat G,\gbel)} \nonumber\\
  & \leq & \left\|K^+_\gbel\right\|_{L^p(\widehat G,\gbel)}\area(\widehat G,\gbel)^{1/q} 
    \leq \cKpo_p \left({\area(\widehat G,\gbel)\over \area(G,\gbel)}\right)^{1/q}.
\end{eqnarray*}
We can apply Proposition~\ref{uupperbound} with $\ka=\cKpo_p$, $r=1/q$, 
$C=\cKpo_p\,\area(G,\gbel)^{-1/q}$ and 
$\mu_0=\mu_1=\area(G,\gbel)$ and we directly get the theorem.
\qed

\section{Cylindrical and conical examples}\label{exsec}

In the previous section we gave a bound on $\osc u$ in terms of 
the $L^p$-norm of the Gaussian curvature $K$, $p>1$, the area and the systole.
Now we will give some examples showing that $\osc u$ is \textbf{not} 
bounded by a function of the $L^1$-norm of $K$, the area and the systole.

In contrast to this, note that the diameter of $\torus$ is bounded by a 
function depending on $\int_{\torus} |K_{\gbel}|$, $\area(\torus,\gbel)$ 
and and $\sys(\torus,\gbel)$, provided that $\int_{\torus} |K_{\gbel}|<4\pi$ 
\cite[Korollar~3.6.8]{ammanndiss}.

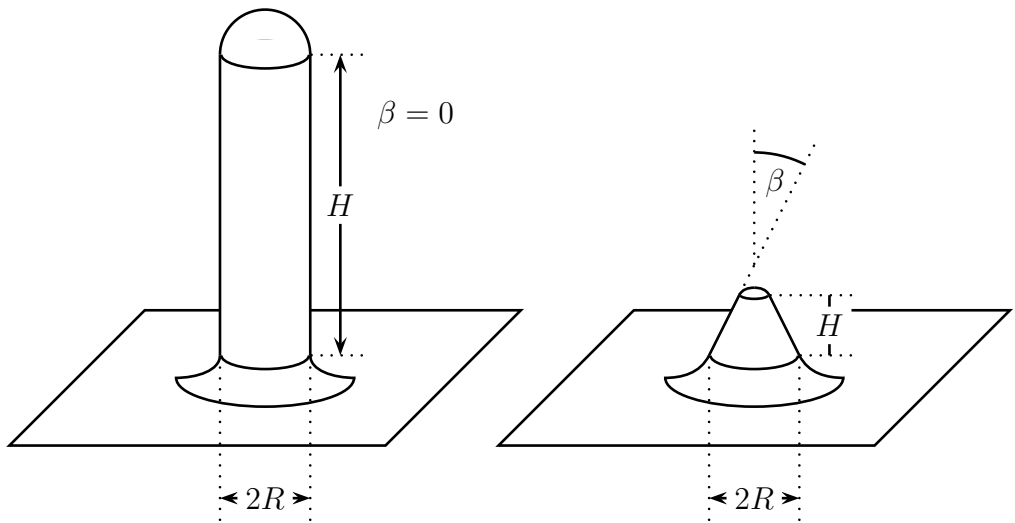
\begin{figure}[htbp]
\def\flaecheneinst{\psset{linewidth=1pt,linecolor=black,linestyle=solid,fillstyle=none}}
\def\querlineinst{\psset{linewidth=1pt,linecolor=black,linestyle=dotted,fillstyle=none}}
\def\pinpfeileinst{\psset{linewidth=1pt,linecolor=black,linestyle=solid,fillstyle=none}}
\def\kegelkappe{
\psbezier(-2,0)(-1.6,0)(-1.2,0.1)(-1,0.5)
\psbezier(2,0)(1.6,0)(1.2,0.1)(1,0.5)
\psline{-}(-1,0.5)(-.333,1.833)
\psline{-}(1,0.5)(.333,1.833)
\psbezier(-.333,1.833)(-.266,1.966)(-.133,2)(0,2)
\psbezier(.333,1.833)(.266,1.966)(.133,2)(0,2)
}

\begin{center}
\psset{unit=1cm}
\begin{pspicture}(-4,-1.7)(11,5.4)

\rput(0,0){
\flaecheneinst
\pspolygon(-3.4,-.9)(1.6,-.9)(3.4,.9)(-1.6,.9)
\psellipse(1.2,0.4)
\psframe*[linecolor=white](-1.2,0)(1.2,.4)
\psbezier(-1.2,0)(-.9,0)(-.6,.1)(-.6,.3)
\psbezier(1.2,0)(.9,0)(.6,.1)(.6,.3)
\psellipse(0,.3)(.6,.2)
\psframe*[linecolor=white](-.6,.3)(.6,4.3)
\psline(-.6,.3)(-.6,4.3)
\psline(.6,.3)(.6,4.3)
\psellipse(0,4.3)(.6,.2)
\psframe*[linecolor=white](-.6,4.3)(.6,4.5)

\querlineinst
\psline(.6,.3)(1.3,.3)
\psline(.6,4.3)(1.3,4.3)
\pinpfeileinst
\psline{<->}(1,.3)(1,4.3)
\rput(1,2.3){\psframebox*[linestyle=none]{$H$}}



\querlineinst
\psline(.6,.3)(.6,-1.9)
\psline(-.6,.3)(-.6,-1.9)
\pinpfeileinst
\psline{<->}(-.6,-1.6)(.6,-1.6)
\rput(0,-1.6){\psframebox*[linestyle=none]{$2R$}}

\flaecheneinst
\psarc(0,4.3){.6}{0}{180}

\rput(2,3.5){$\be=0$}

}

\rput(6.5,0){
\flaecheneinst
\pspolygon(-3.4,-.9)(1.6,-.9)(3.4,.9)(-1.6,.9)
\psframe*[linecolor=white](0,.8)(1.1,1.0)

\psellipse(1.2,0.4)
\psframe*[linecolor=white](-1.2,0)(1.2,.4)
\psellipse(0,.3)(.6,.2)
\psframe*[linecolor=white](-.6,.3)(.6,.5)
\pspolygon*[linecolor=white](-.6,.3)(.6,.3)(0,1.5)
\psellipse(0,1.1)(.2,.066)
\psframe*[linecolor=white](-.2,1.1)(.2,1.3)
\querlineinst
\psline(-.2,1.1)(.8,3.1)
\psline(0,1.5)(0,3.3)
\flaecheneinst
{\psset{unit=.6}\kegelkappe}

\querlineinst
\psline(.6,.3)(.6,-1.9)
\psline(-.6,.3)(-.6,-1.9)
\pinpfeileinst
\psline{<->}(-.6,-1.6)(.6,-1.6)
\rput(0,-1.6){\psframebox*[linestyle=none]{$2R$}}


\rput[bl](.15,2.4){$\be$}
\psarc(0,1.5){1.5}{63.45}{90}

\querlineinst
\psline(.6,.3)(1.3,.3)
\psline(.2,1.1)(1.3,1.1)
\pinpfeileinst
\psline{-}(1,.3)(1,1.1)
\rput(1,.7){\psframebox*[linestyle=none]{$H$}}

}

\end{pspicture}
\end{center}
\caption{cylindrical and conical metric}\label{cynconpic}
\end{figure}

In order to discuss the properties of our examples we will use a lemma 
that can easily be proven by using Lemma~\ref{konformkruem0lemma} 
and the Gauss-Bonnet theorem.

\begin{lemma}
Assume that a disk $D$ carries a rotationally symmetric Riemannian metric $\gbel$
and that in a neighborhood of the boundary $\gbel$ is isometric to a flat ring
of the form $(B_R(0)\ohne B_r(0)\subset \mR^2,\geukl)$.  
Then there is a rotationally symmetric smooth function $u:B_R(0)\to\mR$ vanishing
in a neighborhood of the boundary such that $(D,\gbel)$ is isometric to
$(B_R(0),e^{2u}\,\geukl)$. 
The function $u$ is uniquely determined by these properties. 
\end{lemma}
 
The idea behind the construction of the metrics in this section 
is to start with 
a flat torus $(\torus,\gflach)$, 
to cut out a flat round disk and to replace it 
by a disk $D'$ with a rotationally symmetric metric $g'$. 
Because of the preceeding
lemma the metric $\gbel$ obtained by this replacement 
can be written as $e^{2u}\,\gflach$ where $u$
is a smooth function supported on the disk. Viewed as a function on the disk, 
$u$ is rotationally symmetric. Therefore $\osc u$ can be easily 
estimated using polar coordinates.

The disks $(D',g')$ we glue in are described by figure~\ref{cynconpic}.
For the \emph{cylindrical metric} we construct $(D',g')$ as follows:
we take a cylinder of height $H$ 
and radius $R$, glue it together with a 
half sphere of radius $R$ at one end and a suitable 
socket on the other end. After smoothing we get $(D',g')$.
The resulting metric on $\torus$, the \emph{cylindrical metric}, 
will be denoted $g_{R,H,0}$.  

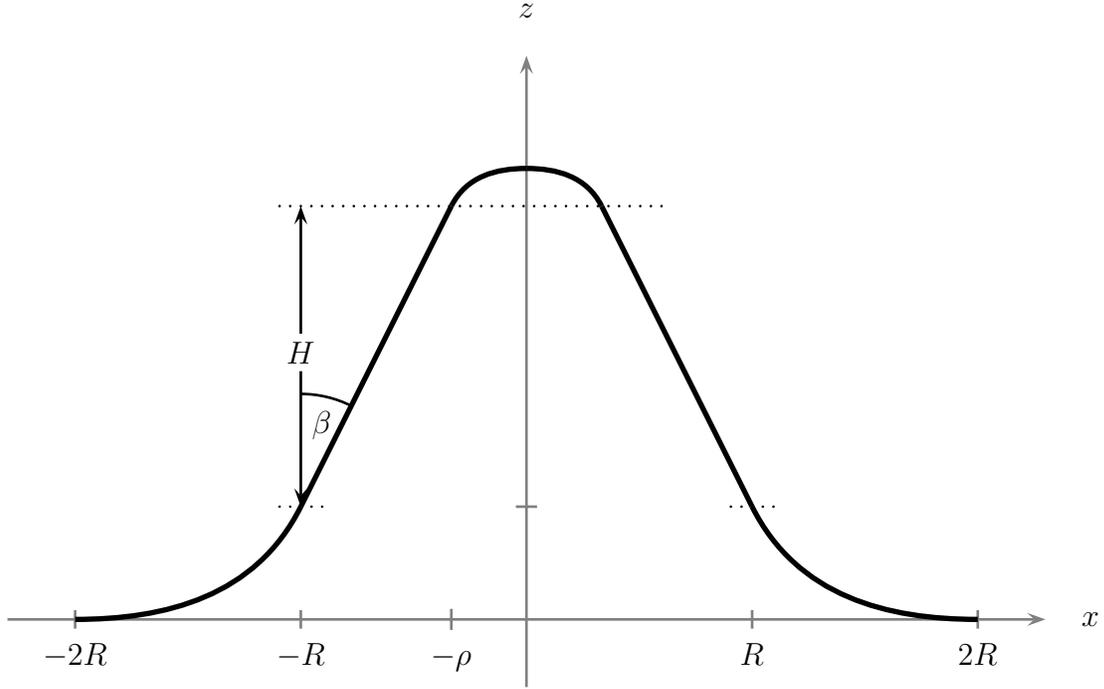
\begin{figure}[tb]
\def\achseneinst{\psset{linewidth=1pt,linecolor=gray,linestyle=solid,fillstyle=none}}
\def\kurveneinst{\psset{linewidth=2pt,linecolor=black,linestyle=solid,fillstyle=none}}
\def\querlineinst{\psset{linewidth=1pt,linecolor=black,linestyle=dotted,fillstyle=none}}
\def\hoehenpfeileinst{\psset{linewidth=1pt,linecolor=black,linestyle=solid,fillstyle=none}}

\def\kegelkappe{
\psbezier(-2,0)(-1.6,0)(-1.2,0.1)(-1,0.5)
\psbezier(2,0)(1.6,0)(1.2,0.1)(1,0.5)
\psline{-}(-1,0.5)(-.333,1.833)
\psline{-}(1,0.5)(.333,1.833)
\psbezier(-.333,1.833)(-.266,1.966)(-.133,2)(0,2)
\psbezier(.333,1.833)(.266,1.966)(.133,2)(0,2)
}

\begin{center}
\psset{unit=3cm}
\begin{pspicture}(-2.5,-0.4)(2.5,3.1)
\achseneinst
\psaxes[labels=none,ticks=x]{->}(0,0)(-2.3,-0.3)(2.3,2.5)
\psline{-}(-.333,4pt)(-.333,-4pt)
\psline{-}(-4pt,.5)(4pt,.5)
\kurveneinst
\kegelkappe

\rput(2.5,0){$x$}
\rput(0,2.7){$z$}
\rput[B](-0.333,-0.2){$-\rho$}
\rput[B](-1,-0.2){$-R$}
\rput[B](1,-0.2){$R$}
\rput[B](-2,-0.2){$-2R$}
\rput[B](2,-0.2){$2R$}
\querlineinst
\psline{-}(-1.1,1.833)(.6,1.833)
\psline{-}(-1.1,0.5)(-0.9,0.5)
\psline{-}(.9,0.5)(1.1,0.5)
\hoehenpfeileinst
\psline{<->}(-1,1.833)(-1,.5)
\rput[b](-1,1.1){\psframebox*[linestyle=none]{$H$}}
\rput[bl](-.95,.8){$\be$}
\psarc(-1,.5){.5}{63.45}{90}
\end{pspicture}
\end{center}
\caption{The conical metric}\label{conicalpic}
\end{figure}

Similarly, for the \emph{conical metric} $g_{R,H,\be}$ 
(Figures~\ref{cynconpic} and
\ref{conicalpic}) we take a 
truncated cone of height $H$, opening angle $\be>0$ 
and the two boundary components are circles of 
radius $R$ and\ $\rho=R-H\sin \be$. 
The end of the truncated cone corresponding to $\rho$ is closed smoothly by a 
topological disk and the other end is put on a socket.

Now we write $g'=e^{2u}\gflach$ and express $u$ in geodesic polar coordinates centered
at the center of $D'$. On the cylinder or cone resp.\ $u$ 
is harmonic and therefore has the form
  $$u(r,\ph)=a + b \log r.$$
The Gauss-Bonnet theorem yields $b=\sin\be -1$.
An elementary calculation shows 
  $$\osc u\geq \left({1\over \sin \be}-1\right)\log\left({R\over \rho}\right)$$
for the conical metric and 
  $$\osc u\geq H/R$$
for the cylindrical metric.

Using Gauss-Bonnet we also get 
  $$\int_{\torus}|K_{g_{R,H,\be}}|=4\pi(1-\sin\be).$$
If $R$ is sufficiently small, the systole does not depend on $H$,
$R$ and $\beta$.

Now for fixed $\be\geq 0$ choose sequences $H_i$ and $R_i$ such that 
$\area(\torus,g_{R_i,H_i,\be})$ is constant and such that $H_i/R_i\to \infty$
or $\rho_i/R_i\to 0$ resp.

So we have constructed families of metrics $g_{R_i,H_i,\be}$ on $\torus$ 
with fixed $\int|K|$, fixed area and fixed systole but $\osc u_i\to \infty$.

\pagebreak

Bernd Ammann\\
Mathematisches Institut\\ 
Universit\"at Freiburg\\
Eckerstr.\ 1\\
79104 Freiburg\\
Germany\\
{\tt ammann@mathematik.uni-freiburg.de}


\begin{thebibliography}{Amm98}

\bibitem[Amm]{ammann:p99b}
B.~Ammann, \emph{Spectral estimates on 2-tori}, Preprint in preparation.

\bibitem[Amm98]{ammanndiss}
B.~Ammann, \emph{{S}pin-{S}trukturen und das {S}pektrum des
  {D}irac-{O}perators}, Ph.D. thesis, University of {F}reiburg, {G}ermany,
  1998, Shaker-Verlag Aachen 1998.

\bibitem[B{\"a}r98]{baer:98}
C.~B{\"a}r, \emph{Extrinsic bounds for eigenvalues of the {D}irac
  operator}, Ann. Global Anal. Geom. \textbf{16} (1998), no.~6, 573--596.

\bibitem[Bes87]{besse:87}
A.~L. Besse, \emph{Einstein manifolds}, Ergebnisse der Mathematik und ihrer
  Grenzgebiete, 3.~Folge, no.~10, Springer-Verlag, 1987.

\bibitem[BM91]{brezis.merle:91}
H.~Brezis and F.~Merle, \emph{Uniform estimates and blow-up behavior for
  solutions of $-{\De} u= {V}(x)\,e^u$ in two dimensions}, Comm. Partial Diff.
  Equat. \textbf{16} (1991), no.~8\ \&\ 9, 1223--1253.

\bibitem[Che81]{chen_bangyen:81}
Bang-Yen Chen, \emph{On the total curvature of immersed manifolds. {V}.
  ${C}$-surfaces in {E}uclidean $m$-space}, Bull. Inst. Math. Acad. Sinica
  \textbf{9} (1981), no.~4, 509--516.

\bibitem[Che98]{chen_xiuxiong:98}
Xiuxiong Chen, \emph{Weak limits of {R}iemannian metrics in surfaces with
  integral curvature bound}, Calc. Var. Partial Differential Equations
  \textbf{6} (1998), no.~3, 189--226.

\bibitem[Gro81]{gromov:81}
M.~Gromov, \emph{Structures m{\'e}triques pour les vari{\'e}t{\'e}s
  {R}iemanniennes}, CEDIC, Paris, 1981.

\bibitem[LS84]{langer.singer:84}
J.~Langer and D.~Singer, \emph{Curves in the hyperbolic plane and mean
  curvature of tori in 3-space}, Bull. London Math. Soc. \textbf{16} (1984),
  531--534.

\bibitem[LY82]{li.yau:82}
P.~Li and S.T. Yau, \emph{A new conformal invariant and its
  applications to the {W}illmore {C}onjecture and the first eigenvalue of
  compact surfaces}, Invent. Math. \textbf{69} (1982), 269--291.

\bibitem[MR85]{montiel.ros:85}
S.~Montiel and A.~Ros, \emph{Minimal immersions of surfaces by
  the first eigenfunctions and conformal area}, Invent. Math. \textbf{83}
  (1985), no.~1, 153--166.

\bibitem[Pin85]{pinkall:85b}
U.~Pinkall, \emph{Hopf tori in ${S}\sp 3$}, Invent. Math. \textbf{81} (1985),
  no.~2, 379--386.

\bibitem[Ros97]{ros:p97}
A.~Ros, \emph{The {W}illmore conjecture in the real projective space},
  Preprint, 1997.

\bibitem[Sim93]{simon:93}
L.~Simon, \emph{Existence of surfaces minimizing the {W}illmore functional},
  Comm. Anal. Geom. \textbf{1} (1993), no.~No.~2, 281--326.

\bibitem[ST70]{shiohama.takagi:70}
K.~Shiohama and R.~Takagi, \emph{A characterization of a standard torus in
  $E^3$}, J. Diff. Geom. \textbf{4} (1970), 477--485.

\bibitem[Tho23]{thomsen:23}
G.~Thomsen, \emph{{\"U}ber konforme {G}eometrie, {I}: {G}rundlagen der
  konformen {F}l{\"a}chentheorie}, Abh. Math. Sem. Hamburg \textbf{3} (1923),
  31--56.

\bibitem[Top98]{topping:p98}
P.~Topping, \emph{Towards the {W}illmore conjecture}, Preprint, 1998.

\bibitem[Wei78]{weiner:78}
J.~L. Weiner, \emph{On a problem of {C}hen, {W}illmore, et al.}, Ind. Math. J.
  \textbf{27} (1978), no.~1, 19--35.

\bibitem[Wil65]{willmore:65}
T.~J. Willmore, \emph{Note on embedded surfaces}, An. \c Sti. Univ. ``Al. I.
  Cuza'' Ia\c si Sec\c t. I a Mat. (N.S.) \textbf{11B} (1965), 493--496.

\bibitem[Wil71]{willmore:71}
T.~J. Willmore, \emph{Mean curvature of {R}iemannian immersions}, J. London
  Math. Soc. (2) \textbf{3} (1971), 307--310.
\end{thebibliography}
\end{document}